\documentclass{article}

\usepackage[numbers]{natbib}
\usepackage{rotating}
\usepackage{graphicx}
\usepackage{amsthm, amssymb, amsmath}

\newtheorem{theorem}{Theorem}

\newtheorem{definition}{Definition}

\newtheorem{Example}{Example}
\newenvironment{example}{\begin{Example}\normalfont}{\end{Example}}

\begin{document}

\title{A Numerical Approach for Solving of Fractional Emden-Fowler Type Equations}

\author{Josef Rebenda
\thanks{CEITEC BUT, Brno University of Technology, Purkynova 123, 612 00 Brno, Czech Republic
(josef.rebenda@ceitec.\allowbreak vutbr.\allowbreak cz).}
\and Zden\v{e}k \v{S}marda
\thanks{CEITEC BUT, Brno University of Technology, Purkynova 123,  612 00 Brno, Czech Republic
(smarda@feec.\allowbreak vutbr.\allowbreak cz).}
}

\date{\copyright\ 2018 AIP Publishing. This article may be downloaded for personal use only. Any other use requires prior permission of the author and AIP Publishing. The following article appeared in "Rebenda, J. and \v{S}marda, Z., A numerical approach for solving of fractional {E}mden-{F}owler type equations, Proceedings of International Conference of Numerical Analysis and Applied Mathematics (ICNAAM 2017), AIP Conference Proceedings, Vol. 1978, 2018" and may be found at https://aip.scitation.org/doi/abs/10.1063/1.5043786}



\maketitle

\begin{abstract}
 In the paper, we utilize the fractional differential transformation (FDT) to solving singular initial value problem of fractional Emden-Fowler type differential equations. The solutions of our model equations are calculated in the form of convergent series with fast computable components. The numerical results show that the approach is correct, accurate and easy to implement when applied to fractional differential equations.
\end{abstract}

\section{INTRODUCTION}
\label{sec1}

Differential equations with fractional order have recently proved to be valuable tools to the modeling of many physical
phenomena \cite{podlubny}-\cite{das}. This is because of the fact that realistic modeling of a physical phenomenon does not depend only on the instant time, but also on the history. This can be successfully achieved by using fractional calculus.

There are many techniques for the solution of fractional differential equations. A good survey of analytical as well as numerical methods is provided in monographs \cite{podlubny}, \cite{miller}, \cite{kilbas}, \cite{diethelm}.

Recently, Adomian decomposition method (ADM) \cite{momani1}-\cite{momani2}, Variational Iteration Method (VIM) \cite{momani1}, \cite{shawa}, Homotopy analysis method (HAM) \cite{wang} belong among the most popular semi-analytical methods.
However, these methods require initial guess or complicated symbolic calculations of  integrals and derivatives. We overcome such drawbacks by implementing simple and easy applicable approach of the fractional differential transformation.

\section{PROBLEM STATEMENT}\label{problemstatement}
In the paper, we apply the fractional differential transformation (FDT) to solving fractional Emden-Fowler type differential
equations in the form
\begin{equation}\label{1}
_0^{C} \! D_t^{2\beta} u + \frac{2}{t^{\beta}}\ {_0^{C}} \! D_t^{\beta} u +f (t)g(u)= 0, \quad t > 0,
\end{equation}
subject to initial conditions
\begin{equation}\label{2}
u(0)=A, \ u'(0) = 0,
\end{equation}
where $\frac{1}{2} < \beta \leq 1$, $A$ is a constant, $f,g$ are continuous functions, $_0^{C} \! D_t^{\lambda}$, $(\lambda >0)$ denotes the fractional derivative of order $\lambda$ in the Caputo sense as defined in the following section. The reason for such special choice of $\beta$ is that the condition $u'(0) = 0$ is used only if $\frac{1}{2} < \beta \leq 1$.

The Emden-Fowler type equations have many applications in the fields of radioactivity cooling and in the mean-field treatment of a phase transition in critical adsorption, kinetics of combustion or reactants concetration in chemical reactor and isothermal gas spheres and thermionic currents \cite{wang}-\cite{smardaz}.

To find a solution of the singular initial value problem for Emden-Fowler type differential equations \eqref{1}, \eqref{2} as well as other various singular initial value problems in quantum mechanics and astrophysics is numerically challenging because of the singular behavior at the origin.

\section{FRACTIONAL DIFFERENTIAL TRANSFORMATION}
In this section, we define the fractional differential transformation (FDT). First we introduce two fractional differential operators.

The fractional derivative in Riemann-Liouville sense is defined by
\begin{equation}\label{3}
_{t_0} \! D^\alpha_{t} f(t) = \frac{1}{\Gamma(n-\alpha)} \frac{d^n}{dt^n}\left[ \int_{t_0}^t \frac{f(s)}{(t-s)^{1+\alpha-n}} ds \right],
\end{equation}
where $n-1 \leq \alpha <n$, $n \in \mathbb{N}$, $t >t_0$. 
 
To avoid fractional initial conditions and to be able to use integer order initial conditions which have a clear physical meaning, we define the fractional derivative in the Caputo sense:
\begin{equation}\label{7}
_{t_0}^{C} \! D_t^{\alpha} f(t) = \frac{1}{\Gamma(n-\alpha)} \int_{t_0}^t \frac{f^{(n)}(s)}{(t-s)^{1+\alpha-n}} ds.
\end{equation}
The relation between the Riemann-Liouville derivative and the Caputo derivative is given by (see e.g. \cite{podlubny}, \cite{kilbas}, \cite{diethelm}) 
\begin{equation}\label{5}
_{t_0}^{C} \! D_t^{\alpha}f(t) =   _{t_0} \! D^\alpha_{t}\left[ f(t) -\sum^{n-1}_{k=0} \frac{1}{k!} (t-t_0)^k f^{(k)}(t_0) \right].
\end{equation}
\begin{definition}\label{def1}
Fractional differential transformation of order $\alpha$ of a real function $u(t)$ at a point $t_0 \in \mathbb R$ in Caputo sense is $^{C} \! \mathcal{D}_{\alpha} \{ u(t) \} [t_0] = \{ U_{\alpha} (k) \}_{k=0}^{\infty}$,
where $k \in \mathbb{N}_0$ and $U_{\alpha} (k)$, the fractional differential transformation of order $\alpha$ of the $(\alpha k)$th derivative of function $u(t)$ at $t_0$, is defined as
\begin{equation}\label{c2}
U_{\alpha} (k) = \frac{1}{\Gamma (\alpha k + 1)} \left[ _{t_0}^{C} \! D_t^{\alpha k} u(t) \right]_{t=t_0},
\end{equation}
provided that the original function $u(t)$ is analytical in some right neighborhood of $t_0$.
\end{definition}
\begin{definition}\label{def2}
Inverse fractional differential transformation of $\{ U_{\alpha} (k) \}_{k=0}^{\infty}$ is defined using a fractional power series as follows:
\begin{equation}\label{c3}
u(t) = {^{C}} \! \mathcal{D}_{\alpha}^{-1} \Bigl\{ \{ U_{\alpha} (k) \}_{k=0}^{\infty} \Bigr\} [t_0]= \sum_{k=0}^{\infty} U_{\alpha} (k) (t-t_0)^{\alpha k}.
\end{equation}
\end{definition}

Convergence of the fractional power series \eqref{c3} in the definition of the inverse FDT was studied in \cite{odibat}. In applications, we will use some basic FDT formulas also listed in \cite{odibat}:
\begin{theorem} \label{t1}
Assume that $\{ F_{\alpha} (k) \}_{k=0}^{\infty}$, $\{ G_{\alpha} (k) \}_{k=0}^{\infty}$ and $\{ H_{\alpha} (k) \}_{k=0}^{\infty}$ are differential transformations of order $\alpha$ of functions $f(t)$, $g(t)$ and
$h(t)$, respectively, and $r, \beta >0$.
\begin{align*}
\text{ If } f(t) &=  (t-t_0)^{r}, \text{ then } F_{\alpha} (k) =\delta \left( k - \frac{r}{\alpha} \right), \text{ where } \delta  
\text{ is the Kronecker delta}. \\
\text{ If } f(t) &= g(t)h(t), \text{ then } F_{\alpha} (k) = \sum_{l=0}^k G_{\alpha} (l) H_{\alpha} (k-l).\\
\text{ If } f(t) &= \frac{g(t)}{(t-t_0)^r}, \text{ then } F_{\alpha} (k) = G_{\alpha} \left( k+ \frac{r}{\alpha} \right).\\
\text{ If } f(t) &= {\displaystyle _{t_0}^{C} \! D_t^{\beta} g(t)}, \text{ then } F_{\alpha} (k) = {\displaystyle \frac{\Gamma(\alpha k + \beta +1)}{\Gamma (\alpha k +1)}} G_{\alpha} \left( k+ \frac{\beta}{\alpha} \right). 
\end{align*}
\end{theorem}


\section{NUMERICAL APPLICATIONS}
Consider singular initial value problem \eqref{1}, \eqref{2}. Applying the FDT, in particular the formulas of Theorem \ref{t1}, to equation \eqref{1}, we obtain the following relation
\begin{equation}\label{e}
\frac{\Gamma(\alpha k+ 2 \beta + 1)}{\Gamma(\alpha k +1)} U_{\alpha} \left(k+ \frac{2 \beta}{\alpha} \right) + 2 \frac{\Gamma(\alpha k+ 2 \beta + 1)}{\Gamma(\alpha k+ \beta + 1)} U_{\alpha} \left(k+\frac{2 \beta}{\alpha} \right)
+ \sum_{l=0}^k F_{\alpha} (l)G_{\alpha} (k-l)=0,
\end{equation}
where $F_{\alpha} (k)$, $G_{\alpha} (k)$ are fractional differential transformations of functions $f(t)$, $g(u)$.

Before we proceed with transformation of initial conditions \eqref{2}, we need to determine the order of the fractional power series $\alpha$. For this purpose, we suppose that $\beta$ is strictly "fractional", i.e. $\beta \in \mathbb{Q}^+$. Then we choose $\alpha$ which satisfies the following conditions:
\begin{enumerate}
 \item
  $0 < \alpha \leq 1$.
 \item
  There is $k_{\beta} \in \mathbb{N}$ such that $\alpha \cdot k_{\beta} = \beta$.
 \item
  There is $k_1 \in \mathbb{N}$ such that $\alpha \cdot k_1 = 1$.
\end{enumerate}
The last condition allows us to use integer order derivatives of $u$ at $t_0$ as initial conditions.

There are infinitely many possibilities for the choice of $\alpha$. However, we propose that $\alpha$ should be chosen as reciprocal of the least common denominator of all orders of fractional derivatives which occur in the considered equation. In our case, we have fractional derivatives of orders $2 \beta$ and $\beta$ in equation \eqref{1}. Recall that we assume $\beta = \frac{p}{q}$ for some $p,q \in \mathbb{N}$. The least common denominator of $\Bigl\{ \frac{2p}{q}, \frac{p}{q} \Bigr\}$ is $q$, and $\alpha = \frac{1}{q}$.

 The transformation of the initial conditions is then defined as
\begin{equation}\label{5'}
U_{\alpha} (k)= \left\{ \begin{array}{ll}
\frac{1}{\Gamma (\alpha k +1)}  \left[ \frac{d^{\alpha k} u(t)}{dt^{\alpha k}} \right]_{t=t_0}, &  {\rm  if} \ \alpha k \in \mathbb{N}, \\
0, & {\rm  if} \ \alpha k \not \in \mathbb{N},
\end{array} \right.
\end{equation}
where $k= 0,1,2,\dots, (\frac{\lambda}{\alpha} -1)$ and $\lambda$ is the order of a considered fractional differential equation, in our case $\lambda = 2 \beta$. In particular, initial conditions \eqref{2} give us $U_{\alpha} (0)=A$ and $U_{\alpha} (q)=0$.
\begin{example}
Consider the following singular initial value problem
\begin{equation}\label{p}
_0^{C} \! D_t^{2\beta} u + \frac{2}{t^{\beta}}\ {_0^{C}} \! D_t^{\beta} u  + u = 0
\end{equation}
subject to intial conditions
$$
u(0)=1,\ u'(0)=0.
$$
We already know that $\alpha =\frac{1}{q}$ and $\beta = \frac{p}{q}$. Then recurrence relation \eqref{e} has the form 
\begin{equation}\label{rec}
\frac{\Gamma(\alpha k+ 2 \beta + 1)}{\Gamma(\alpha k +1)} U_{\alpha} \left(k+ 2p \right) + 2 \frac{\Gamma(\alpha k+ 2 \beta + 1)}{\Gamma(\alpha k+ \beta + 1)} U_{\alpha} \left(k+ 2p \right) +U_{\alpha}(k) = 0.
\end{equation}
From initial conditions we obtain $U_{\alpha} (0)=1$, $U_{\alpha} (1)=0$, \dots, $U_{\alpha} (q-1)=0$, $U_{\alpha} (q)=0$, \dots, $U_{\alpha} (2p-1)=0$. Using the recurrence equation \eqref{rec} we get nonzero coefficients only for $k=0$ and integer multiples of $2p$:
\begin{align*}
k&=0:\ U_{\alpha} (2p)= - U_{\alpha} (0) \left( \Gamma(2 \beta+1)+2\frac{\Gamma(2 \beta+1)}{\Gamma(\beta +1)}\right)^{-1}, \\
k&=1, \dots, 2p-1: \ U_{\alpha} (k+2p)= 0, \\
k&=2p:\ U_{\alpha} (4p) = - U_{\alpha} (2p) \left( \frac{\Gamma(4 \beta+1)}{\Gamma(2 \beta +1)} + 2\frac{\Gamma(4 \beta+1)}{\Gamma(3 \beta +1)} \right)^{-1},\\
k&=2p+1, \ldots, 4p-1 :\ U_{\alpha} (k+2p)=0,\\
\vdots
\end{align*}
Choosing $\beta = 1$ we get the known Lane-Emden type equation
\begin{equation}\label{i}
u'' +\frac{2}{t}u' + u = 0 
\end{equation}
with the exact solution $u(t)= {\displaystyle  \frac{\sin{t}}{t}}$. If we substitute $\alpha = \beta = 1$ in the coefficients $U_{\alpha} (k)$, we have
\begin{align*}
&U_{\alpha} (0) =1,\ U_{\alpha} (1)= 0,\ U_{\alpha} (2) = -\frac{1}{3!},\  U_{\alpha} (3)=0,\  U_{\alpha} (4)  = \frac{1}{5!}, U_{\alpha} (5)=0, \dots,\\
&U_{\alpha} (2k) = \frac{(-1)^{k}}{(2k+1)!}, U_{\alpha} (2k+1) = 0, \dots
\end{align*}
Thus
$$
u(t) = 1 - \frac{t^2}{3!} + \frac{t^4}{5!} - \frac{t^6}{7!} + \dots = \frac{1}{t} \left(t - \frac{t^3}{3!} + \frac{t^5}{5!}- \frac{t^7}{7!} +\dots\right) = \frac{1}{t} \sum\limits_{k=0}^{\infty} \frac{(-1)^{k}}{(2k+1)!} = \frac{\sin{t}}{t}.
$$
We can observe that the solutions of fractional differential equations \eqref{p} converge to the exact solution of differential equation \eqref{i} with the integer order derivative $\beta =1$.
\end{example}


\section{ACKNOWLEDGMENTS}
This research was carried out under the project CEITEC 2020 (LQ1601) with financial support from the Ministry of Education, Youth and Sports of the Czech Republic under the National Sustainability Programme II.


\nocite{*}
\bibliographystyle{aipnum-cp}%
\bibliography{sample}%

\end{document}